\documentstyle[12pt]{article}
\newcommand{\const}{\mathop{\rm const}\limits}

\newcommand{\diam}{\mathop{\rm diam}\limits}

\newcommand{\card}{\mathop{\rm card}\limits}

\newcommand{\Var}{\mathop{\rm Var}\limits}

\newcommand{\radii}{\mathop{\rm radii}\limits}

\begin{document}

\begin{center}

{\bf  UNIFORM MEASURES ON THE  \\

\vspace{4mm}

ARBITRARY COMPACT METRIC SPACES,} \\

\vspace{4mm}

{\bf with applications.}\par

\vspace{3mm}

 $ {\bf E.Ostrovsky^a, \ \ L.Sirota^b } $ \\

\vspace{3mm}

$ ^a $ Corresponding Author. Department of Mathematics and computer science, Bar-Ilan University, 84105, Ramat Gan, Israel.\\

\end{center}

E - mail:\ galo@list.ru; \ eugostrovsky@list.ru\\

\begin{center}

$ ^b $  Department of Mathematics and computer science. Bar-Ilan University,
84105, Ramat Gan, Israel.\\

\end{center}

\vspace{3mm}

E - mail:\ sirota3@bezeqint.net \\

\begin{center}

\vspace{2mm}
                    {\sc Abstract.}\\

 \end{center}

 \vspace{4mm}

  We introduce and investigate in this short report the new notion of uniform measure (distribution)
 on the arbitrary compact metric space.\par
  We consider also some possible applications of these measures in the theory of imbedding theorems and in the theory of
 random processes (fields), in particular, in the so-called  majorizing (and minorizing) measures method,
 belonging to X.Fernique and M.Talagrand.  \par

  These considerations based on the  L.Arnold and P.Imkeller generalization of the  classical
A.M.Garsia-E.Rodemich-H.Jr.Rumsey inequality and X.Fernique-M.Talagrand estimation for random fields.\par

  \vspace{4mm}

 \par

\vspace{4mm}

 \par

\vspace{4mm}

\section{Notations. Definitions. Statement of problem.}

\vspace{3mm}

Let $ (X,d) $ be {\it compact} metric spaces relative a distance (or possible semi-distance)
 $ d = d(x_1, x_2), \ x_1, x_2 \in X. $ The semi-distance $  d  $ differs  on the distance only in that
 $  d(x_1, x_2) = 0 $ does not necessarily imply that $  x_1 = x_2; $  but we can identify all the  points
 $  x_1, x_2 $ for which $  d(x_1, x_2) = 0. $\par

 Let also $ \mu(\cdot) $ be Borelian probability: $  \mu(X) = 1 $ measure on $  X; $ we will consider in this article only
probability Borelian measures. The set of all such a measures will be denoted by $  M\mu = M\mu(X). $ \par

 We denote as usually  $ D = \diam(X) = \max( d(x_1, x_2), \ x_1, x_2 \in X  ) $ the diameter $  X  $ relative the distance $  d; $

 $$
  B(x,r) = B_d(x,r) = \{y: \ y \in X,  \ d(x,y) \le r \}
 $$
 the closed $  d \ - $ ball with center $  x  $ and radii $ r;  \ r \in [0, D]. $\par

  Let us define for any measure $ \mu \in  M\mu(X) $ and for introduced distance  $ d = d(\cdot, \cdot) $
the following two important (for us) functions:

$$
h_-(\mu; r) = h_-(d,\mu; r) = \inf_{x \in X} \mu(B(x,r)),  \eqno(1.0a)
$$

$$
 h_+(\mu; r) = h_+(d,\mu; r) = \sup_{x \in X} \mu(B(x,r)), \  r \in [0, D]. \eqno(1.0b)
$$

 This functions play a very important role in the theory of random processes (fields), see e.g. \cite{Ostrovsky1},  chapter 3,
 p. 95 - 99; \cite{Ostrovsky102}; another applications, in the theory of majorizing measures, will be discuss below. \par

 We define by $ N(T,\epsilon) = N(T, d,\epsilon), \ \epsilon > 0  $ as  usually for the (compact) metric subset $ (T,d), \ T \subset X  $
the minimal number of closed balls with radii $ \epsilon: \  B(x_j, \epsilon)   $ which cover the set $ T: $

$$
T \subset \cup_{j=1}^{N(T)} B(x_j, \epsilon); \ N(\epsilon) := N(X,d,\epsilon) \to \min. \eqno(1.1)
$$

 Recall that the metric space $  (X,d) $ is compact iff it is closed, bounded and $  \forall \epsilon > 0
 \ \Rightarrow N(X,d,\epsilon)  < \infty $ (Hausdorff's theorem). \par

 The (natural) logarithm  of $ N(X,d,\epsilon) $

 $$
H(X,d,\epsilon)  = \log N(X,d,\epsilon)
 $$
is called {\it entropy} of the set $  X $ relative the distance $ d(\cdot, \cdot) $ and also widely used in the entropy
approach to the investigation of continuity of random processes and fields and theory of limit theorems for random processes, see
\cite{Dudley1}, \cite{Fernique1}, \cite{Fernique2}, \cite{Ostrovsky1} etc.\par

 We denote also by $  N(x;\delta, \epsilon) =  N(d,x;\delta, \epsilon)  $ the $  \epsilon \ -  $ entropy of $ \delta \ -  $ ball
with center $  x: $

$$
 N(x;\delta, \epsilon) =  N(d,x;\delta, \epsilon)  = N(B(x,\delta), d, \epsilon);
$$
and correspondingly

$$
 N_-(\delta, \epsilon) =  \inf_{x \in X} N(d,x;\delta, \epsilon), \   N_+(\delta, \epsilon)  = \sup_{x \in X} N(B(x,\delta), d, \epsilon).
$$
 The last notions was introduced by Dmitrovsky \cite{Dmitrovsky1} and was used also in the theory of random fields.\par
 It is clear that

 $$
 \epsilon \ge \delta  \ \Rightarrow N(x;\delta, \epsilon) = 1.
 $$

The similar concept to the entropy forms the so-called capacity $ M = M(X,d,\epsilon) $ of the compact metric space $  X  $ relative to
distance $  d  $  at the point $  \epsilon; \ \epsilon > 0. $  Namely, $  M(X,d,\epsilon) $ is maximal number of {\it disjoint} balls
$ B(y_l, \epsilon) $ belonging to the set $  X; $ and we denote for simplicity

$$
N(\epsilon) = N(X, d, \epsilon), \ \hspace{6mm} M(\epsilon) = M(X, d, \epsilon).
$$

It is known  \cite{Kolmogorov1} that

$$
N(X,d,2 \epsilon) \le M(X,d, \epsilon) \le N(X,d, \epsilon).
$$

 Analogously may be defined the values  $  M(x;\delta, \epsilon) =  M(d,x;\delta, \epsilon)  $ the $  \epsilon \ -  $ capacity
 of $ \delta \ -  $ ball with center $  x: $

$$
 M(x;\delta, \epsilon) =  M(d,x;\delta, \epsilon)  = M(B(x,\delta), d, \epsilon)
$$
and correspondingly

$$
 M_-(\delta, \epsilon) =  \inf_{x \in X} M(d,x;\delta, \epsilon), \   M_+(\delta, \epsilon)  = \sup_{x \in X} M(B(x,\delta), d, \epsilon).
$$

\vspace{3mm}

{\bf  Our purpose is to estimate  both the functions  $ h_{\pm}(\mu; r) $ through the entropy $ H(X,d,\epsilon), \epsilon \in (0, D) $  and
construction an optimal in natural sense  measure, so-called uniform measure, on arbitrary compact metric space. }\par

\vspace{3mm}

 Let $ \epsilon \in (0,D) $ and let $  S(\epsilon) = \{ x_k(\epsilon), \ k = 1,2,\ldots, N(X,d,\epsilon) = N(\epsilon) \} $
 be arbitrary $ d \ - $ optimal $  \epsilon \ - $ net. This net  may be in general case not unique, for example in the case when
 $  X  $ is unit circle on the plane $  R^2 $ equipped with ordinary distance the turn any optimal $ S(\epsilon) $ net on an arbitrary
 angle  remains its optimality.\par
  We can and will conclude by virtue of Egorov's theorem that in the set $  X  $ there is a point $  x_0 $ (the "center" of space)
 for which

 $$
 \sup_{x \in X} d(x, x_0) \le D/2.
 $$

  The value $ D/2 $ may be called {\it radii } of the set (space) $  X:  \radii(X) := D/2. $  Evidently,  under our  agreement

  $$
  \forall \epsilon \ge D/2 \ \Rightarrow H(X,d, \epsilon) = 0.
  $$

 \vspace{3mm}

  Let $ F = \{\nu_{\epsilon}(A) \}, \ \epsilon \in (0,D) $ be a family of discrete Borelian probability measures uniformly
distributed on the set $  S(\epsilon): $

$$
\nu_{\epsilon}(A) = \card  \{x_k(\epsilon) \in S(\epsilon), \ x_k(\epsilon) \in A    \}/N(\epsilon). \eqno(1.2)
$$

 The family $  F  $ satisfies  all the conditions of  famous Prokhorov's theorem \cite{Prokhorov1},
namely, it is weakly compact.  Therefore, there exists a  sequence
$  \nu_{\epsilon_m}(\cdot), \ \epsilon_m > 0, \ \lim_{m \to \infty} \epsilon_m = 0  $  which weakly  converges:

$$
(w)\lim_{ m \to \infty}\nu_{\epsilon_m} = \nu, \eqno(1.3)
$$
where $ \nu $ is also Borelian probability measure on the set $  X. $ \par

 Recall that the equality (1.3) denotes that for arbitrary  continuous bounded function $ f: X \to R $

$$
\lim_{ m \to \infty} \int_X f(x) \ \nu_{\epsilon_m}(dx) = \int_X f(x) \ \nu(dx).
$$

\vspace{3mm}

 {\bf Definition of the uniform measure (distribution) .} \par

\vspace{3mm}

{\it  The arbitrary Borelian probability measure $  \nu $  which is weak limit of some sequence $ \nu_{\epsilon_m} $  as
$ \epsilon_m \to 0+ $  is said to be uniform measure, or equally uniform distribution, defined on the Borelian sets of the space $  X. $  }\par

\vspace{3mm}

 We really proved the existence of uniform  measure (or measures) on the arbitrary metric compact. We investigate further its properties
and consider some applications. \par

\vspace{3mm}

 It is easily to verify  that if $ X $ is compact metrizable topological group, then the uniform measure in our sense is unique and
coincides with the classical normed Haar's measure.  \par

\vspace{4mm}

\section{ Upper estimates. Main result.}

\vspace{3mm}

{\bf Definition 2.1.} {\it The metric space $ (X,d) $ is said to be weakly homogeneous (w.h.), if  there is a constant
$  C_-; \ 0 < C_-  \le 1 $ for which }

$$
\forall \epsilon: 0 < \epsilon \le r \ \Rightarrow  \frac{N_-(r,\epsilon)}{N(\epsilon)}  \ge \frac{C_-}{N(r)}. \eqno(2.1)
$$
 For example, if $  X  $ convex bounded non-trivial domain in the space $ R^l, \  l=1, 2, \ldots  $ with Euclidean distance $ |x-y|  $
satisfying the cone property, and

$$
C_1 |x-y|^{\alpha} \le d(x,y) \le C_2 |x-y|^{\alpha}, \ \alpha = \const \in (0, 1],
$$

$$
 C_1, C_2 = \const, 0 < C_1 \le C_2 < \infty,
$$
where $ |x| $ is ordinary Euclidean norm (or any other equivalent to Euclidean norm),
then the condition (2.1) is satisfied. Moreover,

$$
N(r,\epsilon) \asymp N_-(r,\epsilon) \asymp \left[ \frac{r}{\epsilon} \right]^{l/\alpha}.
$$

\vspace{4mm}

{\bf Theorem 2.1.} {\it Let the space $ (X,d) $ be weakly homogeneous. Then for all the uniform normed Borelian measure } $ \mu  $

$$
h_-(\mu; r) \ge \frac{C}{N(r)}, \ 0 < r  < \infty, \ C = \const \in (0,\infty). \eqno(2.2)
$$

\vspace{3mm}

{\bf Proof.} Let  $ \nu_{\epsilon} $ be discrete measure defined  by the equality (1.2). Then

$$
h_-(\nu_\epsilon,r) = \frac{N_-(r,\epsilon)}{N(\epsilon)}, \ \epsilon \le r.
$$

 It follows from the direct definition of the weak homogeneity that

$$
h_-(\nu_\epsilon,r) \ge \frac{C_-}{N(r)}, \ \epsilon \le r. \eqno(2.3)
$$

Let $ \mu $ be  arbitrary uniform measure on the $  X  $ and
 let $ \epsilon_m, \ m = 1,2,\ldots  $ be any sequence of positive numbers, $ \lim_{m \to \infty} \epsilon_m  = 0 $
such that

$$
\mu = ( w) \ \lim_{m \to \infty} \nu_{\epsilon_m}.
$$

 We deduce passing to the limit as $ \epsilon_m \to 0+ $ in the inequality (2.3)

$$
h_-(\mu,r) \ge \frac{C_-}{N(r)}, \  r > 0; \eqno(2.4)
$$
the case $ r \ge D/2 $   is trivial. \par
 This completes the proof of theorem 2.1.\par

\vspace{4mm}

\section{ Lower estimates. }

\vspace{3mm}

{\bf Proposition 3.1.}  {\it  Let $ \nu(\cdot) $  be additive non - negative probability function defined on the Borelian sets of
the space $  X. $ Then}

$$
\frac{1}{h_+( \nu; 2 \epsilon)} \le M(\epsilon) \le  \frac{1}{h_-(\nu; \epsilon)}. \eqno(3.1)
$$

 {\bf Proof} is alike to one in the coding theory, in the proof of spherical packing boundaries. Indeed,
 let the set $ \{ x_i \}, \ i = 1,2,\ldots, M(\epsilon)  $ be such that

 $$
 X \supset \cup_{i=1}^{M(\epsilon)} B(x_i,\epsilon), \  B(x_i,\epsilon) \cap  B(x_j,\epsilon) = \emptyset, \ i \ne j;
 $$
then

$$
1 = \nu(X) \ge \sum_{i=1}^{M(\epsilon)} \nu(B(x_i, \epsilon)) \ge M(\epsilon) \ h_-(\nu; \epsilon).
$$

 On the other hands, as long as the set $ \{ x_i \}, \ i = 1,2,\ldots, M(\epsilon)  $  is one of the best,

$$
X \subset \cap_{i=1}^{ M(\epsilon)} B(x_i, 2 \epsilon),
$$
 following

 $$
 1 = \nu(X) \le \sum_{i=1}^{ M(\epsilon) } \nu(B(x_i, 2 \epsilon)) \le M(\epsilon) \ h_+(\nu; 2 \epsilon).
 $$

\vspace{4mm}

\section{ Majorizing and  minorizing measures.}

\vspace{3mm}

\begin{center}

{\bf Some applications in the theory of imbedding theorems. }

\end{center}

\vspace{3mm}

 Let $  (Y,\rho) $ be another separable metric spaces, $ m  $ be arbitrary distribution on the Borelian sets $ X, $
 $ f: X \to Y $  be (measurable) function. Let also $ \Phi(z), \ z \in R_+ = [0, \infty) $ be continuous Young-Orlicz function,  i.e.
strictly increasing function such that

$$
\Phi(z) = 0 \ \Leftrightarrow z = 0; \ \lim_{z \to \infty} \Phi(z) = \infty.
$$
 We denote as usually

$$
\Phi^{-1}(w) = \sup \{z, z \ge 0, \  \Phi(z) \le w \}, \ w \ge 0
$$
the inverse function to the function $  \Phi. $

 Let us introduce the Orlicz space $  L(\Phi) = L(\Phi; m \times m, \ X \otimes X)  $ on the set $ X \otimes X $ equipped with
the Young - Orlicz function $ \Phi. $ \par

 \vspace{4mm}

{\it We assume henceforth that for all the values } $ x_1, x_2 \in X, \ x_1 \ne x_2 $ {\it (the case } $ x_1 = x_2 $  {\it is trivial)}
{\it the value  } $ \rho(f(x_1), f(x_2)) $ {\it belongs to the space } $ L(\Phi). $ \par

\vspace{4mm}
 Note that for the existence of such a function $ \Phi(\cdot) $ is necessary and sufficient only the integrability
of the distance   $ \rho(f(x_1), f(x_2)) $ over the product measure $ m \times m:$

$$
\int_X \int_X \rho(f(x_1), f(x_2)) \ m(dx_1) \  m(dx_2) < \infty,
$$
see \cite{Krasnoselsky1}, chapter 2, section 8; where is described a method for building of this function $ \Phi(\cdot). $ \par

 Another natural way, which may gives an optimal up to linear scaling the Orlicz function $  \Phi = \Phi(z). $ Introduce the following functions:

 $$
 \psi(p) := \left[ \int_X \int_X  \rho^p(f(x_1), f(x_2)) \ m(dx_1) \  m(dx_2)   \right]^{1/p}, \ p \ge 1;
 $$

$$
\phi(p) := \left[\frac{p}{\psi(p)} \right]^{-1},
$$

$$
\phi^*(y) := \sup_{p \ge 1} \{p |y| - \phi(p) \} \ -
$$
the Young - Fenchel, or Legendre transform for the function $ \phi(\cdot). $ \par

 The function $  \Phi = \Phi(z) $ may be constructively defined under some natural conditions  as follows:

$$
\Phi(z) := e^{\phi^*(z)} - 1,
$$
see  \cite{Kozachenko1},  \cite{Ostrovsky1}, chapter 1,2. \par

\vspace{3mm}

  Under this assumption the distance $ d = d(x_1, x_2) $ may be constructively defined by the formula:

$$
d_{\Phi}(x_1,x_2) := || \rho(f(x_1), f(x_2))||L(\Phi). \eqno(4.1)
$$

  Since the function $ \Phi = \Phi(z) $ is presumed to be continuous and strictly  increasing, $ \lim_{z \to \infty} \Phi(z) = \infty,  $
it follows from the relation (4.1) that $ V(d_{\Phi}) \le 1, $ where by definition for arbitrary distance function  $ d = d(\cdot, \cdot) $

$$
V = V(d):= \int_X \int_X \Phi \left[ \frac{\rho(f(x_1), f(x_2))}{d(x_1,x_2)} \right] \ m(dx_1) \  m(dx_2). \eqno(4.2)
$$

 Let us define also the following important distance function on the set $  X: $

 $$
   w(x_1, x_2) =  w(x_1, x_2; V)  = w(x_1, x_2; V, m,\Phi) = w(x_1, x_2; V, m,\Phi,d) \stackrel{def}{=}
 $$

$$
 6 \int_0^{d(x_1, x_2)} \left\{ \Phi^{-1} \left[ \frac{4V}{m^2(B(r,x_1))} \right] +
\Phi^{-1}  \left[ \frac{4V}{m^2(B(r,x_2))} \right] \right\} \ dr, \eqno(4.3)
$$
 where $ m(\cdot) $ is probabilistic Borelian measure on the set $ X $  and $ V = V(d). $ \par
  The triangle inequality and other properties of the distance function $ w = w(x_1, x_2) $ are proved in
\cite{Kwapien1}.\par

\vspace{3mm}

{\bf Definition 4.1. } (See  \cite{Kwapien1}). The measure $ m  $ is said to be
{\it minorizing measure } relative the distance $ d = d(x_1,x_2), $ if for each values $ x_1, x_2 \in X
\ V(d) < \infty $ and moreover $ \ w(x_1,x_2; V(d)) < \infty. $\par

\vspace{3mm}

We will denote the set of all minorizing measures on the metric set $ (X,d) $  by  $ \cal{M} = \cal{M}(X).$ \par

 Evidently, if the function $ w(x_1, x_2) $ is bounded, then the minorizing measure $  m  $ is majorizing.  Inverse
proposition is not true, see  \cite{Kwapien1}, \cite{Arnold1}. \par

\vspace{3mm}

{\bf Remark 4.1.} If the measure $  m  $ is minorizing, then

$$
w(x_n, x ; V(d)) \to 0 \ \Leftrightarrow d(x_n, x) \to 0, \ n \to \infty.
$$
 Therefore, the continuity of a function relative the distance $  d  $ is equivalent to
the continuity  of this function  relative the distance $  w.  $ \par

\vspace{3mm}

{\bf Remark 4.2.}  If

$$
\sup_{x_1, x_2 \in X} w(x_1, x_2; V(d)) < \infty,
$$
then the measure $ m $ is called {\it majorizing measure.} This classical definition
with theory explanation and applications basically in the investigation of local structure
of random processes and fields  belongs to
X.Fernique \cite{Fernique1},  \cite{Fernique2},  \cite{Fernique3} and M.Talagrand
\cite{Talagrand1}, \cite{Talagrand2}, \cite{Talagrand3}, \cite{Talagrand4}, \cite{Talagrand5}.
 See also \cite{Bednorz1}, \cite{Bednorz2}, \cite{Bednorz3}, \cite{Dudley1}, \cite{Ledoux1},
 \cite{Ostrovsky100}, \cite{Ostrovsky101}, \cite{Ostrovsky102}. \par

\vspace{4mm}

 The following important inequality belongs to  L.Arnold and P.Imkeller \cite{Arnold1}, \cite{Imkeller1};
see also \cite{Kassman1}, \cite{Barlow1}. \par

\vspace{4mm}

{\bf  Theorem of  L.Arnold and P.Imkeller.} {\it Let the measure $  m  $ be minorizing. Then there exists a
modification of the function $ f $ on the set of zero measure, which we denote also by $ f, $ for which}

$$
\rho(f(x_1), f(x_2)) \le  w(x_1, x_2; V, m,\Phi,d). \eqno(4.4)
$$
{\it  As a consequence: this function $  f  $ is $  d - $ continuous and moreover $ w - $ Lipshitz continuous
with unit constant. }\par

\vspace{4mm}

 The inequality (4.4) of L.Arnold and P.Imkeller is significant generalization of celebrated
Garsia - Rodemich - Rumsey inequality, see \cite{Garsia1}, with at the same applications as mentioned before
\cite{Hu1}, \cite{Ostrovsky100}, \cite{Ostrovsky101}, \cite{Ostrovsky102}, \cite{Ral'chenko1}. \par

 The inequality (4.4) may be also interpreted as an imbedding theorems of the  Sobolev's fractional order space into the space of
continuous functions.\par

\vspace{4mm}

 {\bf   Open question: how to define the majorizing (minorizing) measure (measures)? }\par

 \vspace{4mm}

 Some attempts have been made in the works \cite{Bednorz3}, \cite{Dudley1}, \cite{Fernique3}, \cite{Ledoux1}, \cite{Talagrand2},  \cite{Talagrand4},
 \cite{Talagrand6},  \cite{Talagrand7} etc. \par

 {\it  We intend to offer as the capacity of majorizing measure the (arbitrary)
  uniform distribution on the compact metric space  introduced before. } \par

  By our opinion, this approach is somewhat easier and more convenient than the above.\par

\vspace{3mm}

{\bf Remark 4.3.} The inequality of  L.Arnold and P.Imkeller (4.4) is closely related with the theory of
fractional order Sobolev's - rearrangement invariant spaces, see  \cite{Barlow1}, \cite{Garsia1}, \cite{Hu1}, \cite{Kassman1},
\cite{Nezzaa1}, \cite{Ostrovsky101}, \cite{Ral'chenko1}.\par

\vspace{3mm}

{\bf Remark 4.4.} In the previous articles \cite{Kwapien1}, \cite{Bednorz4} was imposed on the function $ \Phi(\cdot) $
the following $ \Delta^2 $ condition:

$$
\Phi(x) \Phi(y) \le \Phi(K(x+y)), \ \exists K = \const \in (1,\infty), \ x,y \ge 0
$$
or equally

$$
\sup_{x,y > 0} \left[ \frac{\Phi^{-1}(xy)}{\Phi^{-1}(x) + \Phi^{-1}(y)}\right] < \infty. \eqno(4.5)
$$
{\it  We do not suppose this condition. For instance, we can consider the function of a view $ \Phi(z) = |z|^p,  $
which does not satisfy (4.5)}. \par

\vspace{3mm}

{\bf Proposition 4.1.} {\it We deduce using the direct definition of the variable } $ h_-(m,r): \  w(x_1, x_2; V) \le  \overline{w}(x_1, x_2; V),  $
{\it where} $   \overline{w} = \overline{w}(x_1,x_2,V, d,m) = $

$$
\overline{w}(x_1,x_2,V) \stackrel{def}{=} 12 \int_0^{d(x_1, x_2)} \Phi^{-1} \left[ \frac{4V}{h_-^2(m,r)} \right] \ dr.  \eqno(4.6)
$$

\vspace{3mm}

{\bf Proposition 4.2.} {\it We conclude using as the capacity of the measure $  m  $ described below measure $ \nu_{\epsilon}: $ }

$$
6^{-1} w(x_1,x_2,d)  \le \int_{\epsilon}^{d(x_1,x_2)} \Phi^{-1} \left[  \frac{4V}{N^2(\epsilon)} \ N^2(x_2, r,\epsilon)  \right] \ dr+
$$

$$
\int_{\epsilon}^{d(x_1,x_2)} \Phi^{-1} \left[ \frac{4V}{N^2(\epsilon)} \ N^2(x_1, r,\epsilon)  \right] \ dr, \ \epsilon \in (0,D). \eqno(4.7)
$$

\vspace{3mm}

 As a consequence: $ 6^{-1} w(x_1,x_2,d)  \le $

$$
\inf_{\epsilon \in (0,D)} \left\{ \int_{\epsilon}^{d(x_1,x_2)} \left[\Phi^{-1} \left(  \frac{4V}{N^2(\epsilon)} \ N^2(x_2, r,\epsilon)  \right) +
 \Phi^{-1} \left(  \frac{4V}{N^2(\epsilon)} \ N^2(x_1, r,\epsilon) \right) \right] \ dr \right\}. \eqno(4.8a)
$$

or $ 6^{-1} w(x_1,x_2,d)  \le $

$$
\underline{\lim}_{\epsilon \to 0+} \left\{ \int_{\epsilon}^{d(x_1,x_2)} \left[\Phi^{-1} \left(  \frac{4V}{N^2(\epsilon)} \ N^2(x_2, r,\epsilon)  \right) +
 \Phi^{-1} \left(  \frac{4V}{N^2(\epsilon)} \ N^2(x_1, r,\epsilon) \right) \right] \ dr \right\}. \eqno(4.8b)
$$

\vspace{3mm}

{\bf Proposition 4.3.} {\it Suppose in addition that the (compact) metric space $ (X,d)  $ is weak homogeneous. If we choose as a capacity
of the measure  $ m $ arbitrary uniform distribution on the set $  X, $ then  }

$$
w(x_1,x_2,V) \le 12 \int_0^{d(x_1, x_2)} \Phi^{-1} \left[ \frac{4V N^2(r)}{C_-^2} \right] \ dr.  \eqno(4.9)
$$

\vspace{3mm}

\section{ Application to the problem of boundedness and continuity of random fields. }

\vspace{3mm}

\begin{center}

{\bf General Orlicz and majorizing measures method approach.}\\

 \end{center}

\vspace{4mm}

  Let $ \xi = \xi(x), \ x \in X $ be separable numerical stochastic continuous  (i.e. continuous in probability  relative the distance $  d ) $
  random field (r.f.),  defined aside the probability space  on the introduced before set $ X, $  not necessary to be Gaussian.\par

   The correspondent  set of all elementary events,
   probability and expectation we will denote as ordinary correspondingly by $ \Omega, \ {\bf P}, \ {\bf E,}  $  and  the probabilistic
 Lebesgue-Riesz $ L_p $   norm of a random variable (r.v) $  \eta $  we will denote as follows:

 $$
 |\eta|_p \stackrel{def}{=} \left[ {\bf E} |\eta|^p \right]^{1/p}.
 $$

   We find in this section some sufficient condition for boundedness of $ \xi(x) $ (a.e.).

    Recall that the first publication about fractional  Sobolev's  inequalities
\cite{Garsia1}  was devoted in particular to the such a problem; see also articles \cite{Hu1}, \cite{Ostrovsky101}, \cite{Ral'chenko1}. \par

\vspace{3mm}

 Let $ \Phi = \Phi(u) $ be again the Young-Orlicz function. We impose together with M.Ledoux and M.Talagrand,
  \cite{Ledoux1}, chapter 11 the following conditions on the function $ \Phi(\cdot): $

 $$
 \Phi^{-1}(xy) \le C \left(  \Phi^{-1}(x) + \Phi^{-1}(y)  \right), \ \exists C = \const < \infty
 $$

 $$
 \int_0^1 \Phi^{-1} (x^{-1}) \ dx < \infty. \eqno(5.0)
 $$

 We will denote the Orlicz norm by means of the function $  \Phi $
 of a r.v. $ \kappa $ defined on our probabilistic space $ (\Omega, {\bf P}) $ as $ |||\kappa|||L(\Omega,\Phi) $   or for simplicity
 $ |||\kappa|||\Phi. $  \par
  We introduce  the so-called {\it natural} (new) distance $ d_{\Phi}(x_1,x_2) $ on the set $   X  $ as follows:

  $$
   d_{\Phi} = d_{\Phi}(x_1,x_2):= |||\rho(\xi(x_1),\xi(x_2))|||L(\Omega,\Phi), \ x_1,x_2 \in X.  \eqno(5.1)
  $$

\vspace{4mm}

 {\bf Definition 5.1., see \cite{Ledoux1}, chapter 11.} \par

\vspace{3mm}

 The Borelian probabilistic measure $  m(\cdot) $ on the Borelian  sets $  X  $  equipped by the distance $  d_{\Phi} = d_{\Phi}(\cdot, \cdot)  $
 is called {\it  majorizing, } iff

$$
\gamma_m(X, d, \Phi) := \sup_{x \in X} \int_0^D \Phi^{-1} \left( \frac{1}{m(B_{d_{\Phi}}(x,r)) } \right) \ dr < \infty. \eqno(5.2)
$$

\vspace{4mm}

{\bf Theorem 5.1., see \cite{Ledoux1}, section 11.2.} {\it  Let  $ m(\cdot) $ be the probabilistic  majorizing measure on the set $  X  $
relative the distance $  d_{\Phi}(\cdot, \cdot). $  Then }

$$
 \sup_{y \in X} \left\{ {\bf E} \sup_{x \in X} [ \xi(x) - \xi(y)] \right\} \le K(\Phi) \ \gamma_m(X, d_{\Phi}, \Phi). \eqno(5.3)
$$

  {\it As  a consequence:  the random field   $  \xi = \xi(x), \ x \in X  $ is bounded with probability one. } \par

\vspace{3mm}

{\bf Remark 5.1.} M.Ledoux and M.Talagrand proved in addition that if

$$
\lim_{\delta \to 0+} \sup_{x \in X} \int_0^{\delta} \Phi^{-1} \left( \frac{1}{m(B_{d_{\Phi}}(x,r)) } \right) \ dr =0,
$$
then the r.f. $ \xi = \xi(x) $ is $ d_{\Phi} \ - $ continuous with probability one.\par

\vspace{4mm}

  M.Ledoux and M.Talagrand used in particular the following  majorizing measure:

 $$
  m_{N} := \sum_{l> l_0} 2^{ -l + l_0 }  N \left( X,d_{\Phi}, 2^{-l}\right) \sum_{x \in S(2^{-l}) }  \ \delta_x,
 $$
where $ \delta_x $ denotes an usual unit Dirac measure concentrated at the point $ x. $ \par
 This gives the following entropic estimate:

$$
\sup_{y \in X} \left\{ {\bf E} \sup_{x \in X} [ \xi(x) - \xi(y)] \right\} \le K_2(\Phi) \
\int_0^D \Phi^{-1} [ N(X,d_{\Phi},r) ] \ dr . \eqno(5.4)
$$

\vspace{3mm}

{\bf Definition 5.2.} The r.f. $ \xi = \xi(x) $ defined on the compact metric set $  X = (X,d)  $
is said to be {\it  quasi - homogeneous } relative some probability Borelian measure $  m, $ if

 $$
 \sup_{r > 0} \left[ h_+(m,r) /h_-(m,r) \right] < \infty.  \eqno(5.5)
 $$

 R.M.Dudley in \cite{Dudley1}, p. 59-62 proved that if  the r.f. $  \xi = \xi(x) $ is
quasi - homogeneous relative some probability Borelian measure $  m, $ then

$$
\int_0^D \Phi^{-1} [ N(X,d_{\Phi},r) ] \ dr  < \infty \ \Leftrightarrow \gamma_m(X, d_{\Phi}, \Phi) < \infty.
$$

 M.Talagrand in \cite{Talagrand3} proved that if the r.f. $ \xi = \xi(x) $ is Gaussian and centered, then the finiteness of  the value
$ \gamma_m(X, d_{\Phi}, \Phi) $ for some Borelian probability measure $  m  $ on the set $  X  $ is necessary and sufficient condition
for it boundedness. \par
 Here evidently
 $$
  \Phi(z) = \Phi_2(z) \stackrel{def}{=} \exp(z^2/2) - 1.
 $$
Therefore, if in addition  the Gaussian r.f. $ \xi = \xi(x) $ is quasi - homogeneous, the condition

$$
\int_0^D \Phi^{-1} [ N(X,d_{\Phi_2},r) ] \ dr  < \infty  \eqno(5.6)
$$
 is necessary and sufficient  for a.e. boundedness of $ \xi(\cdot) $ and moreover for it continuity relative the
natural Dudley's distance

$$
d_D(x_1,x_2) = \left[ \Var(\xi(x_1) - \xi(x_2)  \right]^{1/2},
$$
which is linear equivalent of course in the considered here Gaussian case to the distance $ d_{\Phi_2}. $

 \vspace{3mm}

 We intent to improve the estimate (5.4) by means of using uniform distribution on the compact space $ (X, d_{\Phi}), $
as well as  in the fourth section.  \par

\vspace{3mm}

{\bf Proposition 5.1.} {\it We deduce using again the direct definition of the variable } $ h_-(d,m,r):   $

$$
\gamma_m(X, d_{\Phi}, \Phi) \le  \int_0^D \Phi^{-1} \left( \frac{1}{h_-(d_{\Phi}, m, r) } \right) \ dr.  \eqno(5.7)
$$

\vspace{3mm}

{\bf Proposition 5.2.} {\it We conclude using as the capacity of the measure $  m  $ described below measure $ \nu_{\epsilon}: $ }

$$
\gamma_m(X, d_{\Phi}, \Phi)  \le \int_{\epsilon}^{D} \Phi^{-1} \left( \frac{ N(d_{\Phi}, r,\epsilon)}{N(d_{\Phi},\epsilon)}  \right) \ dr, \
 \epsilon \in (0,D). \eqno(5.8)
$$

\vspace{3mm}

 As a consequence:

$$
\gamma_m(X, d_{\Phi}, \Phi)  \le \inf_{\epsilon \in (0,D)}
\int_{\epsilon}^{D} \Phi^{-1} \left( \frac{ N(d_{\Phi}, r,\epsilon)}{N(d_{\Phi}, \epsilon)}  \right) \ dr, \
 \epsilon \in (0,D); \eqno(5.8a)
$$

\vspace{3mm}

$$
\gamma_m(X, d_{\Phi}, \Phi)  \le \underline{\lim}_{\epsilon \to 0+}  \int_{\epsilon}^{D}
\Phi^{-1} \left( \frac{ N(d_{\Phi}, r,\epsilon)}{N(d_{\Phi},\epsilon)} \right) \ dr, \
 \epsilon \in (0,D). \eqno(5.8b)
$$

\vspace{3mm}

{\bf Proposition 5.3.} {\it Suppose in addition that the (compact) metric space $ (X,d)  $ is weak homogeneous. If we choose as a capacity
of the measure  $ m $ arbitrary uniform distribution on the set $  X, $ then  }

$$
\gamma_m(X, d_{\Phi}, \Phi) \le \int_0^{D} \Phi^{-1} \left( \frac{N(d_{\Phi}, r)}{C_-} \right) \ dr.  \eqno(5.9)
$$

\vspace{3mm}

 Let us show the other approach,  which nay gives us the {\it exponential estimates} for distribution of $ \sup_{ x \in X} | \xi(x) |, $
see e.g.   \cite{Ostrovsky102}, \cite{Ostrovsky105}. \par

  Let $ \Phi = \Phi(u) $ be again the Young-Orlicz function, generated as before by means of a {\it family } of random variables
 $ \rho(\xi(x_1),\xi(x_2)), \ x_1, x_2 \in X.  $  We will denote the Orlicz norm by means of the function $  \Phi $
 of a r.v. $ \kappa $ defined on our probabilistic space $ (\Omega, {\bf P}) $ as $ |||\kappa|||L(\Omega,\Phi) $   or for simplicity
 $ |||\kappa|||\Phi. $  \par

  We can and will suppose without loss of generality

 $$
\sup_{x_1, x_2 \in X} |||  \rho(\xi(x_1),\xi(x_2)) |||\Phi  = 1.
 $$

  We introduce also a new {\it non-random} so-called  {\it natural}  bounded: $ v_{\Phi}(x_1,x_2) \le 1   $
 distance on the set $  X, \  v = v_{\Phi}(x_1,x_2), $ i.e. generated only by the random field $ \{\xi(x) \}, \ x \in X  $ as follows:

  $$
  v := v_{\Phi} = v_{\Phi}(x_1,x_2):= |||\rho(\xi(x_1),\xi(x_2))|||L(\Omega,\Phi), \ x_1,x_2 \in X.  \eqno(5.10)
  $$

\vspace{4mm}

{\bf Theorem 5.1.} {\it  Let  $ m(\cdot) $ be the probabilistic minorizing measure on the set $  X  $
relative the distance $  d_{\Phi}(\cdot, \cdot). $ There exists a non-negative random variable $  Z = Z(v_{\Phi},m) $
with unit expectation: $ {\bf E} Z = 1 $  for which }

$$
\rho(\xi(x_1), \xi(x_2)) \le  \overline{w}(x_1, x_2; Z(v_{\Phi},m)). \eqno(5.11)
$$
  {\it As  a consequence:  the r.f.  $  \xi = \xi(x)  $ is $ d \ - $ continuous with probability one. } \par
\vspace{4mm}

{\bf Proof.}  We pick

$$
Z = \int_X \int_X \Phi \left( \frac{\rho(\xi(x_1), \xi(x_2))}{v_{\Phi}(x_1,x_2)}  \right) \ m(dx_1) \ m(dx_2).
$$
We have by means of theorem Fatou - Tonelli

$$
 {\bf E} Z  = {\bf E} \int_X \int_X  \Phi \left( \frac{\rho(\xi(x_1),\xi(x_2))}{v_{\Phi}(x_1,x_2)} \right) \ m(dx_1) \ m(dx_2) =
$$

$$
 \int_X \int_X {\bf E} \Phi \left( \frac{\rho(\xi(x_1),\xi(x_2))}{v_{\Phi}(x_1,x_2)} \right) \ m(dx_1) \ m(dx_2) = 1, \eqno(5.12)
$$
since $ \int_X \int_X  m(dx_1) \  m(dx_2) = 1.  $\par
 It remains to use our proposition 4.1 and  apply the L.Arnold and P.Imkeller inequality. \par

\vspace{4mm}

\section{Non-asymptotical exact exponential estimates for distribution of maximum of random fields.}

\vspace{3mm}

 \begin{center}

 {\bf Grand Lebesgue spaces approach.}

 \end{center}

\vspace{3mm}

 Let  $ \xi = \xi(x), \ x \in X  $ be again separable random field (or process)
with values in the real axis $  R, \ T = \{x \} $ be arbitrary Borelian subset of  $ X.$

 Denote

 $$
 \overline{\xi}_T = \sup_{x \in T} \xi(x), \hspace{6mm}  \overline{\xi} = \overline{\xi}_X = \sup_{x \in X} \xi(x), \eqno(6.0)
 $$

$$
Q(u) = Q_{\xi}(u) = Q(X,u) \stackrel{def}{=} {\bf P}( \sup_{x \in X} \xi(x) > u), \ u \ge 2; \eqno(6.1)
$$

$$
Q_+(u) = Q_{\xi, +}(u) = Q_+(X,u) \stackrel{def}{=} {\bf P}( \sup_{x \in X} |\xi(x)| > u), \ u \ge 2; \eqno(6.1a)
$$
 and for arbitrary Borelian subset $ T \subset X $ we denote

$$
Q(T,u) = {\bf P}( \sup_{t \in T} \xi(t) > u), \ u \ge 2. \eqno(6.2)
$$

$$
Q_+(T,u) = {\bf P}( \sup_{t \in T} |\xi(t)| > u), \ u \ge 2. \eqno(6.2a)
$$

 Of course,  here

$$
Y = R, \ \rho(y_1,y_2) = |y_1 - y_2|.
$$

  Our purpose in this section is obtaining an {\it exponentially exact} as $ u \to \infty, \ u > u_0 = \const > 0 $
 estimation for the probability $ Q(u) \stackrel{def}{=} Q(X,u),  $ in the modern
 terms of "minorizing measures"  and the so-called $ B(\phi) $ spaces, which are a particular case of Grand Lebesgue  spaces.  \par
 We intent to improve and simplify foregoing results using method of majorizing measures.\par

  In the entropy terms this problem is considered in  \cite{Dudley1}, \cite{Fernique2}, \cite{Fernique3},
 \cite{Ostrovsky1}, chapter 3,  \cite{Ostrovsky103}: in the more modern terms of majorizing measures- in \cite{Ledoux1},
\cite{Ostrovsky102},  \cite{Talagrand1}  etc. Note that the method of majorizing measures, or equally generic chaining method,
gives more strong results  as entropy technique
in the theory of continuity of random processes and fields, see \cite{Fernique1} - \cite{Fernique3},
\cite{Talagrand1} - \cite{Talagrand7}, \cite{Ostrovsky102} etc. \par

  The estimations of $  Q(u) \stackrel{def}{=} Q(X,u)  $ are used in the Monte-Carlo method, statistics,
 numerical methods etc., see \cite{Frolov1}, \cite{Grigorjeva1},
 \cite{Ostrovsky1},  \cite{Ostrovsky105}, \cite{Ostrovsky106}. \par

\vspace{3mm}

{\bf A. Preliminaries.} \\

{\it We  suppose in this section that the random field $ \xi(x)  $ to be centered and satisfies the uniform Kramer's condition,
so that the natural function, equally the uniform (exponential) moment generating function

$$
\phi(\lambda) = \log \sup_{x \in X} {\bf E} \exp(\lambda \ \xi(x)) \eqno(6.3)
$$
is finite in some non-trivial interval } $ \lambda \in (-\lambda_0, \lambda_0), \ \lambda_0 = \const \in (0, \infty].  $\par
 Then we may introduce the following Young-Orlicz function (up to multiplicative positive constant)

 $$
 \Phi_{\phi}(u) = \exp(\phi^*(u))-1,
 $$
so that  $ \sup_{x \in X} ||\xi(x)||B(\phi) = 1  $ and following  $ \sup_x ||\xi(x)||(\Phi) < \infty. $ \par
 We recall first of all some propositions from the theory  of the so-called $  B(\phi) \ $  spaces; more detail
descriptions see, e.g. in  \cite{Kozachenko1}, \cite{Ostrovsky1}, \cite{Ostrovsky102}.\par

 Let $ \phi = \phi(\lambda), \lambda \in (-\lambda_0, \lambda_0), \ \lambda_0 = \const \in (0, \infty] $ be some even strong
convex which takes positive values for positive arguments twice continuous differentiable function, such that
$$
 \phi(0) = 0, \ \phi^{//}(0) \in(0,\infty), \ \lim_{\lambda \to \lambda_0} \phi(\lambda)/\lambda = \infty.
$$

 We say that the {\it centered} random variable (r.v) $ \xi = \xi(\omega) $
belongs to the space $ B(\phi), $ if there exists some non-negative constant
$ \tau \ge 0 $ such that

$$
\forall \lambda \in (-\lambda_0, \lambda_0) \ \Rightarrow
{\bf E} \exp(\lambda \xi) \le \exp[ \phi(\lambda \ \tau) ]. \eqno(6.4)
$$
 The minimal value $ \tau $ satisfying (6.4) is called a $ B(\phi) \ $ norm
of the variable $ \xi, $ write
 $$
 ||\xi||B(\phi) = \inf \{ \tau, \ \tau > 0: \ \forall \lambda \ \Rightarrow
 {\bf E}\exp(\lambda \xi) \le \exp(\phi(\lambda \ \tau)) \}.
 $$
 This spaces are very convenient for the investigation of the r.v. having a
exponential decreasing tail of distribution, for instance, for investigation of the limit theorem,
the exponential bounds of distribution for sums of random variables,
non-asymptotical properties, problem of continuous of random fields,
study of Central Limit Theorem in the Banach space etc.\par

  The space $ B(\phi) $ with respect to the norm $ || \cdot ||B(\phi) $ and
ordinary operations is a Banach space which is isomorphic to the subspace
consisted on all the centered variables of Orlicz's space $ (\Omega,F,{\bf P}), N(\cdot) $ with $ N \ - $ function

$$
N(u) = \exp(\phi^*(u)) - 1, \ \phi^*(u) = \sup_{\lambda} (\lambda u -\phi(\lambda)).
$$
 The transform $ \phi \to \phi^* $ is called Young-Fenchel transform. The proof of considered assertion used the
properties of saddle-point method and theorem of Fenchel-Moraux: $ \phi^{**} = \phi. $

 The next facts about the $ B(\phi) $ spaces are proved in \cite{Kozachenko1}, \cite{Ostrovsky1}, p. 19-40:

$$
{\bf 1.} \ \xi \in B(\phi) \Leftrightarrow {\bf E } \xi = 0, \ {\bf and} \ \exists C = \const > 0,
$$

$$
U(\xi,x) \le \exp(-\phi^*(Cx)), x \ge 0,
$$
where $ U(\xi,x)$ denotes in this article the {\it tail} of
distribution of the r.v. $ \xi: $

$$
U(\xi,x) = \max \left( {\bf P}(\xi > x), \ {\bf P}(\xi < - x) \right),
\ x \ge 0,
$$
and this estimation is in general case asymptotically exact. \par

\vspace{3mm}

{\bf 2.}  The $  B(\phi) $ norm  $ ||\eta||B\phi  $ of a random variable $ \eta $
is linear equivalent to the following norm
$$
 ||\eta||G\psi \stackrel{def}{=} \sup_{p \ge 1} \left[  \frac{|\eta|_p}{\psi(p) }  \right]
$$
in the Grand Lebesgue space $ G\psi: $

$$
C_1 ||\eta||G\psi   \le ||\eta||B\phi  \le C_2 ||\eta||G\psi,
$$
where

$$
\psi(p) = \frac{p}{\phi^{-1}(p)}.
$$

\vspace{3mm}

{\bf 3.}  For the non-centered random variable $  \kappa $ with exponentially decreasing tail of distribution the $ B(\phi) $ norm
may be defined, for instance, as follows:

$$
||\tau||B(\phi) := \left[ ||\tau - {\bf E} \tau||^2  + ( {\bf E} \tau)^2   \right]^{1/2}.
$$

\vspace{3mm}

{\bf B. Auxiliary constructions. } \\

\vspace{3mm}

Let the function $ \overline{w} = \overline{w}(x_1, x_2, Z) $ was  defined by equality (4.6).
 Let  $ m(\cdot) $ be the probabilistic minorizing measure on the set $  X  $
relative the distance $  d_{\Phi}(\cdot, \cdot), $ for example, uniform measure. We know, see theorem (5.1) that
there exists the non-negative random variable $  Z  $ with unit expectation: $ {\bf E} Z = 1 $  for which

$$
\rho(\xi(x_1), \xi(x_2)) \le  \overline{w}(x_1, x_2; Z).
$$

 Define the random variable (or equally function of $  Z ) \ \theta = \theta(Z) $ as follows:

$$
\theta = \gamma(Z) := \sup_{x_1, x_2 \in X} \overline{w}(x_1, x_2; Z),
$$
and  introduce  the following semi-distance function

$$
\zeta(x_1,x_2) := \sup_{ Y > 0 } \left[ \frac {\overline{w}(x_1, x_2; Y)}{\gamma(Y)} \right],
$$
so that

$$
\overline{w}(x_1, x_2; Y) \le \zeta(x_1,x_2) \cdot \gamma(Y)
$$
or equally

$$
\overline{w}(x_1, x_2; Z) \le \zeta(x_1,x_2) \cdot \theta. \eqno(6.5)
$$

\vspace{4mm}

{\bf C. Lemma 6.1.}  {\it Assume in addition that the condition (6.3) is satisfied.  Then we have in the representation }
(6.5)

$$
||\theta||B(\phi) = K < \infty. \eqno(6.6)
$$

 {\bf Proof} follows immediately from the main result of paper \cite{Ostrovsky107}, (see also the article \cite{Ostrovsky102},)
 where using for us result it is formulated in the similar language of $ G\psi $ spaces. \par

\vspace{4mm}

{\bf D. Some estimates.}\\

\vspace{3mm}

 Let $  B_{\zeta}(x_0,\delta) $ be a $ \zeta \ -  $ ball on the set $  X  $ with the center at the point $  x_0 $  and radii $ \delta  >: 0: $

$$
B_{\zeta}(x_0,\delta) = \{ x, \ x \in X, \ \zeta(x,x_0) \le \delta \}.
$$
 The diameter of the set $  X  $ relative the distance $  \zeta $ will be denoted by  $ D_{\zeta}; $ so that we can restrict the
values $ \delta $  inside the interval $ \delta \in (0, D_{\zeta}). $ \par

 We derive using triangle inequality  for the values $ x $ inside the ball $ B_{\zeta}(x_0, \delta):  $

 $$
 ||\xi(x)||B(\phi) \le ||\xi(x_0) ||B(\phi) + K \cdot \delta \le 1 + K \ \delta,
 $$
therefore

$$
\sup_{x_0 \in X} Q_{B_{\zeta}(x_0,\delta)}(u) \le \exp \left( - \phi^*\left( \frac{u}{1 + K \ \delta}   \right)   \right). \eqno(6.7)
$$

\vspace{4mm}

{\bf  E. Concluding exponential estimation.}\\

\vspace{3mm}

  Denote by $  S_{\zeta}(\epsilon) $ the minimal $ \epsilon \ - $ net for the set $  X  $ relative the distance $  d_{\zeta}: $

$$
S_{\zeta}(\epsilon)  = \{x_1(\epsilon), x_2(\epsilon), \ldots, x_{N(\epsilon)}(\epsilon) \};
$$
 where  by our notations (temporarily, in this section)

 $$
 \card \left(S_{\zeta}(\epsilon) \right) = N(X,\zeta,\epsilon) \stackrel{def}{=} N(\epsilon).
 $$
 We  have using (6.7):

 $$
Q(X,u) = {\bf P}\left( \sup_{x \in X}  \xi(x) > u \right) =
 {\bf P} \left( \cup_{j = 1}^{ N(\delta)} \{ \sup_{x \in B(x_j(\delta))}  \xi(x) > u  \}  \right) \le
 $$

$$
\sum_{j=1}^{N(\delta)}  {\bf P} \left( \sup_{x \in B(x_j(\delta))}  \xi(x) > u   \right) =
\sum_{j=1}^{N(\delta)} Q_{B_{\zeta}(x_j(\delta),\delta)}(u) \le
$$

$$
N(X,\zeta,\delta) \cdot  \exp \left( - \phi^*\left( \frac{u}{1 + K \ \delta}   \right)   \right). \eqno(6.8)
$$

\vspace{3mm}

{\bf E. Total:}\\

$$
Q(X,u) \le
\inf_{\delta \in (0, D_{\zeta} ) } \left\{ N(X,\zeta,\delta) \cdot  \exp \left( - \phi^*\left( \frac{u}{1 + K \ \delta}   \right) \right) \right\}, \eqno(6.9)
$$
and correspondingly

$$
Q_+(X,u) \le 2
\inf_{\delta \in (0, D_{\zeta} ) } \left\{ N(X,\zeta,\delta) \cdot  \exp \left( - \phi^*\left( \frac{u}{1 + K \ \delta} \right)  \right) \right\}, \eqno(6.9a)
$$

\vspace{4mm}

\begin{center}

{\bf Examples.}

\end{center}

\vspace{3mm}

{\bf Example 6.1.} Suppose in addition  to the conditions of proposition 4.2 analogously to the work \cite{Ostrovsky102}
that the function $ u \to \Phi(u), \ u > 0 $ is logarithmical convex:

$$
(\log \Phi)^{''} (u) > 0.
$$
 Let also $ \gamma = \const \in (0,1), \ C_1 = \const \in (0,\infty).  $ Denote

 $$
 \delta_0 = \delta_0(u; \gamma, \Phi) = \frac{C_1 \gamma}{u \cdot [\log \Phi]'(u)}. \eqno(6.10)
 $$
 We obtain choosing $ \delta = \delta_0 $ substituting into (6.9) that for all sufficiently greatest values
 $ u: \ \delta_0(u; \gamma, \Phi) < D_{\zeta} $

$$
Q(X,u) \le \frac{(1-\gamma)^{-1}}{\Phi(u)} \cdot N \left( \frac{C_2 \gamma}{u \cdot  [\log \Phi]'(u)} \right). \eqno(6.11)
$$

\vspace{4mm}

{\bf Example 6.2.} Let now $ \Phi(u) = \exp(u^2/2) - 1 $ (subgaussian case).  Suppose

$$
N(\epsilon) \le C_3 \epsilon^{-\kappa}, \ \epsilon \in (0,D), \ \kappa = \const > 0.  \eqno(6.12)
$$
 The value $ \kappa $ is said to be {\it majorital} dimension of the set $ X $ relative the distance $ w. $ \par
The optimal value $ \gamma $ in (4.20) if equal to $ \gamma = \gamma_0 := \kappa/(\kappa + 1) $ and we conclude for the
values $ u $ such that

$$
\delta_0 = \frac{C_2 \kappa}{(\kappa + 1) u^2} \le D_{\zeta}:
$$

$$
Q(X,u) \le C_3 \ C_2^{-\kappa} \ \kappa^{-\kappa} \ (\kappa+1)^{ \kappa + 1} \  u^{2 \kappa} \ e^{-u^2/2}. \eqno(6.13)
$$

\vspace{4mm}

{\bf Example 6.3.} Assume that  in the example 6.2 instead the condition (6.12) the following condition holds:

$$
N(\epsilon) \asymp  C_4 e^{ C_5 \ \epsilon^{-\beta}}, \ \epsilon \in (0,D);  \ \beta  = \const > 0.  \eqno(6.14)
$$

 Then

 $$
 Q(X,u) \le e^{-0.5 u^2 + C_6 u^{2\beta/(\beta+1)} }, \ u \ge C_7. \eqno(6.15)
 $$

 Note that in the  case  $ \beta \ge 2 $ the so-called entropy series

$$
\sum_{n=1}^{\infty} 2^{-n} \ H^{1/2} \left(X,w, 2^{-n} \right)
$$
diverges.\par

\vspace{4mm}

\section{Concluding remarks.}

\vspace{3mm}

{\bf A. Degrees.} \par
 Let $ X = [0,1]^n,  \ n = 2,3,\ldots. $  In the articles   \cite{Ral'chenko1}, \cite{Hu1}  is obtained a multivariate generalization of
famous Garsia-Rodemich-Rumsey inequality \cite{Garsia1}. Roughly speaking, instead degree "2" in the inequalities (4.2) and
(4.3) stands degree 1  and coefficients dependent on $ d. $ \par
 The ultimate value of this degree in general case of arbitrary metric space $ (X,d) $ is now unknown; see also  \cite{Arnold1},
\cite{Imkeller1}. \par

\vspace{3mm}

{\bf B. Spaces.} \par
 Notice that in all considered cases and under our conditions when
  $ \sup_x ||\xi(x)||B(\phi) < \infty,  $ then $ ||\sup_x \xi(x)||B(\phi)  < \infty. $ But in the article
\cite{Ostrovsky104} was constructed a "counterexample": there exists a continuous a.e random
process for which

$$
 \sup_x ||\xi(x)||B(\phi)  < \infty, \hspace{5mm}  ||\sup_x \xi(x)||B(\phi)  = \infty.
$$
 This circumstance imply that our conditions are  only sufficient but not necessary.\par

\vspace{3mm}

{\bf C.  Lower bounds.} \par
 The lower estimates for probabilities $ Q(u)  $ are obtained e.g. in \cite{Ostrovsky1}, chapter 3, sections 3.5-3.8.
They are  obtained in entropy terms, all the more so in the terms of minorizing measures.  \par
Note that the lower bounds in \cite{Ostrovsky1} may coincide up to multiplicative constants with upper bounds. \par

\vspace{3mm}

{\bf D. Other possible applications.} \par

We  can obtain having in hand the exact  exponential  estimations for tail  function for maximum distribution of random field,
the sufficient conditions for weak compactness of random fields in the space of continuous functions, and in particular
derive the sufficient conditions for Central Limit Theorem in Banach space, alike ones in the article \cite{Ostrovsky102}. \\

\end{document}